\documentclass[11pt,a4paper,english]{article}
\usepackage[a4paper,left=3cm,right=2cm]{geometry}
\usepackage{babel}
\usepackage[latin1]{inputenc}
\usepackage{amsmath}
\usepackage{amsthm}
\usepackage{amsfonts}
\usepackage{indentfirst}
\usepackage{graphicx}
\newtheorem{de}{Definition}
\newtheorem{pro}{Proposition}
\newtheorem{teo}{Theorem}
\newtheorem{rem}{Remark}
\newtheorem{lem}{Lemma}
\newtheorem{exa}{Example}

\newcommand{\co}{{\mathcal O}}

\newcommand{\gp}{\mathbb{P}}

\newcommand{\gz}{\mathbb{Z}}

\newcommand{\gc}{\mathbb{C}}

\newcommand{\cp}{{\mathcal P}}

\newcommand{\findemo}{$\ \ \square$}
\newcommand{\cf}{{\mathcal F}}

\renewcommand{\int}{{\rm int}}

\newcommand{\Pp}{{\mathbb P}}

\title{ {\bf On the Characterization of Algebraically Integrable Plane Foliations}}
\author{C. Galindo \thanks{Supported by the Spain Ministry of Education
 MTM2007-64704 and Bancaixa P1-1A2005-08 }
 \and F.  Monserrat$^{*}$}
\date{}

\begin{document}
\maketitle

\section{Introduction}
At the end of the 19th century, Darboux \cite{dar}, Poincar\'e
\cite{poi1,poi2,poi3}, Painlev\'e \cite{pai} and Autonne \cite{aut}
studied the problem that, nowadays, can be stated as follows: to
characterize those algebraic foliations on the projective plane over
the field of complex numbers (plane foliations, in the sequel) which
have a rational first integral. Most recent contributions related to
that problem are focused to solve the so called Poincaré problem,
that consists of bounding the degrees of the irreducible algebraic
leaves of a foliation in terms of numerical data related to it
\cite{car, ca-ca, l-n, es-kl}. From this bound, one can try to get a
rational first integral through an algebraic computation.


In this note, we only consider non-degenerated plane foliations
${\cal F}$. When ${\cal F}$ has a rational first integral $R$, the
singular points of ${\cal F}$ are, on the one hand, the
indeterminacies of the rational map $f_R: \gp^2\cdots \rightarrow
\gp^1$ which $R$ defines ($Ind(R)$) and, on the other hand, the
singular points of the fibers of this rational map ($Sing(R)$). In
the case considered in this paper, non-degenerated foliations,
$Ind(R)$ ($Sing(R)$, respectively) coincides with the set of {\it
non-reduced} ({\it reduced}, respectively) singularities of ${\cal
F}$, denoted by $NRed({\cal F})$ ($Red({\cal F})$, respectively).

In Proposition \ref{primero}, for a plane foliation ${\cal F}$ of
degree different from 1, we consider the set of rational functions
$R$ such that $Ind(R)=NRed({\cal F})$ and $Sing(R)=Red({\cal F})$
and characterize when an element of that set  is a rational first
integral of ${\cal F}$ in terms of local conditions on the
singularities of the closures of the fibers of $f_R$. These
conditions only depend on the analytic type of the mentioned
singularities. An important ingredient of the proof is the fact,
proved by Campillo and Olivares in \cite{cam-oli}, that foliations
of degree different from one are determined by the scheme associated
with the indeterminacies of the polarity map of the foliation.



If we know the Seidenberg's resolution of ${\cal F}$ \cite{sei}, the
ideas used to prove the above result can be strengthened to get our
main result, Theorem \ref{principal}, which is a characterization
theorem for foliations of degree different from 1 having a rational
first integral. The characterization conditions are  two numerical
ones (related to  the eigenvalues of the non-reduced singularities)
and other conditions that have to do with certain sheaf on the
surface obtained after Seidenberg's resolution process. Note that
when ${\cal F}$ is of degree one and all its singularities  are
non-degenerated, with rational quotient of eigenvalues, then it has
always a rational first integral.

In the previous results, an important object is the set $NRed({\cal
F})$. For this reason, we end this paper by proving in Theorem
\ref{cota} that, under the assumption that ${\cal F}$ is
algebraically integrable, the cardinality of $NRed({\cal F})$ is
larger than the degree of ${\cal F}$. Notice that this means that
the degree of a non-degenerated foliation, $r$, provides the minimum
number, $r+1$, of points in $\gp^2$ through which infinitely many
leaves of the foliation go.

\section{Plane foliations and algebraic integrability}

Along this note, $\gp^2$ will denote the complex projective plane,
$\co_{\gp^2}$ ($\co_{\gp^2}^{an}$, respectively) its structural
sheaf as algebraic variety over $\gc$ (its sheaf of holomorphic
functions, respectively). A {\it plane foliation} or a {\it
foliation} on $\Pp^2$, $\mathcal{F}$, will be an (algebraic)
foliation (with singularities and singular set of codimension two)
of degree $r$ ($r \geq 0$) on $\Pp^2$. That is, $\mathcal{F}$ is a
nontrivial map of vector bundles $ \mathcal{F}: H^{\otimes(-r+1)}
\rightarrow T \Pp^2$ whose cokernel is torsion-free, where $T
\Pp^2$ and $H$ are the corresponding bundles to the tangent and
hyperplane sheaves on $\Pp^2$ and $H^{\otimes n}$ means the
$n$-fold tensor product of $H$ if $n >0$, the $(-n)$-fold tensor
product of $H^{\vee}$ if $n <0$ and the trivial line bundle if
$n=0$.

A plane foliation $\mathcal{F}$  provides a tangent direction (or,
equivalently, a projective line) to all points in $\Pp^2$ but
finitely many called {\it singularities of} $\mathcal{F}$, set
denoted by $ Sing(\mathcal{F})$. This allows to define the so-called
{\it polarity map} $\phi: \Pp^2 \cdots \rightarrow (\Pp^2)^\vee$,
that is, the rational map which sends each point $p \in \Pp^2
\setminus Sing (\mathcal{F}) $ to that point in $(\Pp^2)^\vee$
corresponding to the projective line that $\mathcal{F}$ associates
with $p$. The scheme-theoretic fibers $\phi^* E$ of the lines $E$ in
$(\Pp^2)^\vee$ are degree $-(r+1)$ curves on $\Pp^2$ and so, since
lines in $(\Pp^2)^\vee$ can be identified with points in $\Pp^2$, we
can assign to each point $p \in \Pp^2$ a curve of degree $r+1$
called the {\it polar of $p$ with respect to $\mathcal{F}$}. The
indeterminacy ideal $\mathcal{J}$ of the polarity map gives a
subscheme of $\gp^2$, named the {\it singular subscheme} of
$\mathcal{F}$, which determines the foliation \cite{cam-oli}.

$\mathcal{F}$, in analytic terms, can be defined, up to
multiplication by a non-zero complex number, by a reduced
homogeneous vector field (up to a suitable multiple of the radial
vector field) or, equivalently, by a reduced homogeneous 1-form,
$\Omega \equiv A \; dX + B \; dY + C \; dZ$, $(X:Y:Z)$ being
projective coordinates on $\Pp^2$, such that $A,B,C$ are homogeneous
polynomials of degree $(r+1)$ in $X,Y,Z$ without common factors and
satisfying the Euler condition $X A + Y B +Z C =0$. This last
condition is what allows to recover $\Omega$ from local data, which
will be considered over the corresponding affine charts of $\Pp^2$.
For example, $\Omega$ over the chart $Z \neq 0$ will be expressed by
$\omega = a(x,y) dx+ b(x,y) dy$, where $x= X/Z$, $y=Y/Z$, $a(x,y)
=A(X,Y,1)$, $ b(x,y)= B(X,Y,1)$. As a consequence a point $p$ in the
above chart will be a singularity of $\mathcal{F}$ if, and only if,
$a(p) =b(p) =0$. Recall that the colength of the stalk
$\mathcal{J}_p$ of the above defined indeterminacy ideal (which is
the ideal $(a, b){\mathcal O}_{\Pp^2, p}$ in the coordinates
$\{x,y\}$) in the local ring ${\mathcal O}_{\Pp^2, p}$, $\mu_p({\cal
F})$, is called the {\it Milnor number} of $\mathcal{F}$ at the
point $p$.

Fixed homogeneous coordinates  on $\gp^2$, a rational function $R$
of $\gp^2$ is defined by the quotient of two homogeneous polynomials
in these coordinates without common  components and with the same
degree, $F$ and $G$. Associated with $R$, we consider the linear
pencil $\cp(R)$ consisting on the projective curves with equations
$\alpha F+\beta G=0$, where $(\alpha, \beta)$ runs over
$\gc^2\setminus \{0\}$. The curves in $\cp(R)$ are the topological
closures of the fibers of the rational map $\gp^2\cdots \rightarrow
\gp^1$ that $R$ defines. The base locus of $\cp(R)$ is the set of
indeterminacies of this map. As we have said, it will be denoted by
$Ind(R)$. Also, we shall denote by $Sing(R)$ the set of singular
points of the curves of the pencil $\cp(R)$ which are not in
$Ind(R)$.


By definition, a foliation $\cf$ has a {\it rational (or
meromorphic) first integral} if there exists a rational function $S$
of $\gp^2$ such that $dS \wedge \Omega=0$. In this case, there
exists a rational function $R=F/G$ such that all the rational first
integrals of ${\cal F}$ have the form $\frac{P(F,G)}{Q(F,G)}$, where
$P,Q$ are homogeneous polynomials of the same degree; such a
function $R$ will be called a {\it primitive rational first
integral} and it satisfies that its associated pencil of plane
curves $\cp(R)$ is irreducible (that means that it has irreducible
general elements). Moreover, the leaves of ${\cal F}$ are the
algebraic curves whose irreducible components are components of some
curve of the mentioned pencil.

We only consider foliations $\cf$, which we shall call {\it
non-degenerated} ones, such that the Milnor numbers of all its
singularities equal one. For each singularity $p$ of $\cf$, the pair
$(\delta,\rho)$ of eigenvalues of the linear part of any local
vector field defining ${\cal F}$ at $p$ is named the (pair of) {\it
eigenvalues associated with $p$}. Notice that $\delta \neq 0 \neq
\rho$. We are interested in algebraic integrability; so we assume
that $\delta$ and $\rho$ are integers and that
$\gcd(\delta,\rho)=1$, after dividing by its greatest common divisor
if it is necessary. When $\delta/\rho$ is not positive then we shall
say that ${\cal F}$ has a {\it reduced} singularity at $p$. We shall
denote by $Red({\cal F})$ ($NRed({\cal F})$, respectively) the set
of reduced (non-reduced, respectively) singularities.

A well known fact, that we shall use systematically along the paper,
is that if ${\cal F}$ is a non-degenerated foliation with rational
first integral $R$, then $Ind(R)=NRed({\cal F})$ and
$Sing(R)=Red({\cal F})$.

\section{Characterization of foliations with a rational first integral}

In this section, we  study the problem of characterizing when a
non-degenerated plane foliation has a rational first integral. First
of all, we recall some definitions that we shall use in the sequel.

\begin{de}
{\rm A polynomial $P\in \gc[x,y]$ is  {\it weighted-homogeneous}
with weights $(w_1,w_2)$, where $w_1,w_2$ are fixed rational
numbers, if it can be expressed as a linear combination of monomials
$x^{\alpha}y^{\beta}$ for which
$\frac{\alpha}{w_1}+\frac{\beta}{w_2}=1$. }
\end{de}

\begin{de}
{\rm A plane curve singularity is {\it quasi-homogeneous} if it is
analytically equivalent to a singularity whose local equation is
defined by a weighted-homogeneous polynomial. }
\end{de}

Stand $\mathbb{Z}_+$ for the set of positive integers. Given $a,b,k
\in \mathbb{Z}_+$ such that $\gcd(a,b)=1$, we shall denote by
$S(a,b,k)$ the topological singularity type of a plane curve
singularity with equation
$x^{ka}+y^{kb}=0$, where $\{x,y\}$ are local coordinates.
We shall also consider the following set of topological singularity
types:
$${\cal S}=\{S(a,b,k)\mid a,b,k\in \gz_+ \mbox{ and } \gcd(a,b)=1\}.$$

The topology of a quasi-homogeneous singularity determines the
weights of all quasi-ho\-mo\-ge\-neous polynomials defining the
singularity \cite{yosu}. Therefore, the quasi-homogeneous
singularities with singularity type $S(a,b,k)$ correspond to the
weights $(ka,kb)$.

\begin{de}
{\rm We shall say that a plane curve singularity  is {\it nodal} if
its local equations have the form $g_1^n g_2^m=0$, where $m,n$ are
positive integers and $g_1,g_2$ define analytically irreducible
regular  germs which are transversal. }
\end{de}

\begin{de}
{\rm Let ${\cal F}$ be a non-degenerated foliation on $\gp^2$ and
$p$ a non-reduced singular point of ${\cal F}$ with associated
eigenvalues $(\delta, \rho)$. ${\cal F}$ is said to be {\it
linearizable} at $p$ if there exist analytical coordinates $\{u,v\}$
at $p$ such that ${\cal F}_p$ is defined by the differential form
$\rho v du-\delta u dv$.

}
\end{de}

Although the results of the next lemma are well-known for the
specialists, we include its proof for the reader's convenience.

\begin{lem}\label{lema}
Let ${\cal F}$ be a non-degenerated foliation on $\gp^2$ admitting a
rational first integral $R$. Then,
\begin{itemize}
\item[(a) ] ${\cal F}$ is linearizable at every non-reduced
singularity.

\item[(b) ] A point $p\in \gp^2$ is a reduced singularity of
${\cal F}$ if, and only if, it is not a base point of $\cp(R)$ and
the unique curve $H$ of $\cp(R)$ passing through $p$ has a nodal
singularity at this point.

\end{itemize}
\end{lem}

\noindent {\it Proof.} For non-resonant singularities, $(a)$ follows
from Poincaré's linearization Theorem (see, for instance,
\cite{cs1}) and otherwise, bearing in mind that the singularity is
non-reduced, from  the Poincaré-Dulac normal form Theorem (see
\cite{cs1} and \cite[pages 16 and 17]{brun}).


To prove $(b)$, assume that $p$ is a reduced singularity. Set $h:= u
\prod_{i=1}^s h_i^{e_i}=0$ a local equation of the germ at $p$ of
$H$,
where $u$ is a unit and $h_i=0$  an analytically irreducible germ
for $i=1,2,\ldots,s$. Then, $d(u^{-1}h)= w \prod_{i=1}^s
h_i^{e_i-1}$,
where $w=adx+bdy$  with $a,b\in \co_{\gp^2,p}^{an}$ and
$\gcd(a,b)=1$,  $\{x,y\}$ being local coordinates (see \cite[Section
2.4]{julio}). Therefore, the multiplicity of $\cf$ at $p$ (i.e., the
minimum of the $m_p$-adic orders of $a$ and $b$, $m_p$ being the
maximal ideal of $\co_{\gp^2,p}^{an}$) is $\nu=-1 + \sum_{i=1}^s
\nu(h_i)$, where $\nu(h_i)$ denotes the multiplicity of the germ
$h_i$. Since $\nu=1$, one has that $\sum_{i=1}^s \nu(h_i)=2$ and, in
particular, either $s=1$ or $s=2$.

If $s=1$, then $h$ is $u(x^2+f)^{e_1}$ for some local coordinates
$\{x,y\}$,  $f\in (x,y)^3$ and $u$  a unit. This shows a local
expression of $\cf$ at $p$, what, clearly, cannot hold since ${\cal
F}$ is non-degenerated. If $s=2$, either $h=0$ defines a nodal
singularity, or $h=u(x+f_1)^{e_1}(x+f_2)^{e_2}$  with $f_1,f_2\in
(x,y)^2$, which again contradicts the fact that $\cal F$ is
non-degenerated.




The converse is clear since, in certain coordinates $\{x,y\}$, a
local equation for $H$ is $h=0$, with $h=(x+f)^{e_1}(y+g)^{e_2}$ and
$f,g\in (x,y)^2\subseteq \co_{\gp^2,p}^{an}$, and  hence, $p$ is a
reduced singularity of ${\cal F}$. \findemo

\medskip

Recall that if $C$ is a curve, $p$ a point in $C$, $\eta:
\overline{C} \rightarrow C$ the normalization map and $\gamma$ the
class map, that is the composition on the canonical sheaf $\gamma:
\Omega_C^1 \rightarrow \eta_* \Omega_{\overline{C}}^1 \rightarrow
\omega_C$, the {\it Milnor} ({\it Tjurina}, respectively) {\it
number} of $C$ at $p$ is $\mu(C,p):=l({\rm Cok}(\gamma \circ d)_p)$
($\tau(C,p): = l({\rm Ext}^1_{{\cal O}_{C,p}} (\Omega_{C,p}^1, {\cal
O}_{C,p}))$, respectively), where $d: {\cal O}_{C} \rightarrow
\Omega_C^1$ is the universal derivation and $l$ the length.

\begin{pro} \label{primero}
Let ${\cal F}$ be a non-degenerated plane algebraic foliation of
degree different from 1 and  let $R$ be a rational function such
that $Ind(R)=NRed({\cal F})$ and $Sing(R)=Red({\cal F})$. Then, $R$
is a rational first integral of ${\cal F}$ if and only if the
following conditions are satisfied:
\begin{itemize}

\item[1. ] The singularity at every point in $Sing(R)$ of the
curve of $\cp(R)$ passing through it is nodal.

\item[2. ] There exist two curves of the pencil $\cp(R)$ (say,
with equations $H_1=0$ and $H_2=0$) such that, for each $p\in
Ind(R)$, the germs at $p$ of the curves with equations $H_1=0$ and
$H_2=0$ are equisingular and the one of $H_1H_2=0$ satisfies the
following properties:

\begin{itemize}

\item[(a) ] it is a reduced germ whose topological singularity
type belongs to ${\cal S}$,

\item[(b) ] its associated Milnor and Tjurina numbers coincide.

\end{itemize}

\end{itemize}

\end{pro}

\noindent {\it Proof}. Assume that $R$ is a rational first integral
for ${\cal F}$. Condition 1 is satisfied by Part $(b)$ of Lemma
\ref{lema}.

In order to prove 2, pick two distinct general elements (with
equations $H_i=0$, $i=1,2$) of the pencil $\cp(R)$) such that the
local equation at each point $p\in Ind(R)$ of the curve $C$ defined
by $H_1H_2=0$ is reduced. Consider a point $p\in Ind(R)$. By Part
$(a)$ of Lemma \ref{lema}, there exist local coordinates $u,v\in
\co_{\gp^2,p}^{an}$ such that the foliation ${\cal F}$ is defined
locally at $p$ by the differential 1-form $\rho vdu-\delta u dv$,
where $(\delta, \rho)$ is the pair of eigenvalues associated with
$p$. If $h_1$ ($h_2$, respectively) in $\co_{\gp^2,p}^{an}$ defines
the germ  provided by $H_1$ ($H_2$, respectively), since $H_1 = 0$
and $H_2 = 0$ are general elements of $\cp(R)$ one has that $h_1$
and $h_2$ give equisingular reduced germs (and, obviously, without
common factors). Moreover, since $\frac{u^{\rho}}{v^{\delta}}$ is a
primitive local rational first integral of the germ of foliation
${\cal F}_p$, there exist a unit $z$ and homogeneous polynomials in
two variables $P$ and $Q$ of the same degree and without common
factors such that $h_1=zP(u^{\rho},v^{\delta})$ and
$h_2=zQ(u^{\rho},v^{\delta})$ (see, for instance, \cite[Section
2.9]{julio}).  $z^{-2}h_1h_2$ factorizes into a product of
polynomials of the type $L(u^{\rho},v^{\delta})$, $L$ being an
homogeneous polynomial of degree $1$. Hence, it is clear that
Condition $2(a)$ is satisfied. Finally, notice that $h_1h_2$ defines
a quasi-homogeneous singularity and, then, Condition $2(b)$ follows
from the second ``Satz'' in \cite[page 123]{saito}.

Conversely, consider the foliation ${\cal H}$ provided by the
derivation of the rational function $R$. Taking into account the
hypotheses of the statement, the equality $Sing({\cal H})=Sing({\cal
F})$ clearly holds. Next, we shall prove the equality ${\cal
F}={\cal H}$ and, hence, the result. By applying \cite[Th.
3.5]{cam-oli}, it is enough to show that ${\cal H}$ is a
non-degenerated foliation.

Let $H_1$ and $H_2$ be as in Condition 2 and take whichever point
$p\in Ind(R)$. Consider the curve $C$ of equation $H_1 H_2=0$. By
Condition 2(b) and the second ``Satz'' in \cite[page 123]{saito},
$C$ has a quasi-homogeneous singularity at $p$. Therefore, there
exist local analytic coordinates $\{u,v\}$ and a quasi-homogeneous
polynomial $P$ such that $P(u,v)=0$ defines this singularity.
Moreover, applying Condition 2(a), its singularity type belongs to
${\cal S}$ and there exist three positive integers $\rho, \delta$
and $k$ such that $\gcd(\rho,\delta)=1$ and $(k\delta, k\rho)$ is
the pair of weights corresponding to $P(u,v)$. It is straightforward
to see that $P(u,v)=Q(u^{\rho},v^{\delta})$, where $Q$ is an
homogeneous polynomial in two variables of degree $k$. Therefore,
$P(u,v)=\prod_{i=1}^{k} L_i(u^{\rho},v^{\delta})$, $L_i$ being
distinct homogeneous polynomials of degree 1. Now, $k$ is even and
$P(u,v)=Q_1(u^{\rho},v^{\delta})Q_2(u^{\rho},v^{\delta})$, where
$Q_1$ and $Q_2$ are homogeneous polynomials of degree $k/2$ and
$Q_1(u^{\rho},v^{\delta})=0$ ($Q_2(u^{\rho},v^{\delta})=0$,
respectively) is an equation of the germ at $p$ of the curve $H_1=0$
($H_2=0$, respectively) in the local coordinates $u$ and $v$ (notice
that these facts are true because the mentioned germs are
equisingular). Taking into account that the derivation of
$Q_1(u^{\rho},v^{\delta})/Q_2(u^{\rho},v^{\delta})$ defines the same
local foliation as the germ ${\cal H}_p$, we deduce that
$u^{\rho}/v^{\delta}$ is a primitive local first integral of ${\cal
H}_p$ and, therefore, $\mu_p({\cal H})=1$. Finally, Condition 1
ensures that $\mu_p({\cal H})=1$ for all $p\in
Sing(R)$, which concludes the proof.\findemo\\


Let ${\cal F}$ be a non-degenerated plane foliation.
By Seidenberg's resolution process \cite{sei}, there exists a
sequence of point blowing-ups whose composition, $\pi: Y \rightarrow
\gp^2$, satisfies that the foliation on $Y$ $\pi^*({\cal F})$ has
only reduced singularities (see \cite[pages 12 and 13]{brun}). Such
a composition morphism is called a {\it minimal resolution} of $\cal
F$ if it is minimal with respect to the number of involved
blowing-ups. We shall denote by $\pi_{\cal F}:X_{\cal F}\rightarrow
\gp^2$ a minimal resolution of singularities of ${\cal F}$ and by
${\cal C}_{\cal F}$ the set of centers of the blowing-ups that are
involved in it (notice that this set is, essentially, unique).
${\cal C}_{\cal F}$ is a disjoint union $\bigcup {\cal C}_{{\cal
F},p}$, $p$ running over $NRed({\cal F})$, where ${\cal C}_{{\cal
F},p}:=\{q\in {\cal C}_{\cal F} \mid q\geq p\}$ and $\geq$ is the
partial ordering on ${\cal C}_{\cal F}$  given by $q\geq r$ if and
only if $q$ is infinitely near to $r$ \cite[3.3]{casas}.

Now we shall define a family of sheaves on $X_{\cal F}$ associated
with this minimal resolution. For each $p\in NRed({\cal F})$ we
distinguish two cases:
\begin{itemize}
\item[-] ${\cal F}$ is linearizable at $p$. Then there exist local
coordinates $u,v\in \co_{\gp^2,p}^{an}$ such that ${\cal F}_p$ is
defined by the differential 1-form $\rho v du-\delta u dv$,
$(\delta, \rho)$ being the pair of eigenvalues associated with $p$.
Let $J_p$ be the ideal of $\co_{\gp^2,p}^{an}$ generated by
$u^{\rho}$ and $v^{\delta}$. Elementary computations show that the
infinitely near points needed to eliminate the base points of the
ideal $J_p$ (see \cite[7.2]{casas}) coincide with those in ${\cal
C}_{{\cal F},p}$. As a consequence, ${\cal C}_{{\cal F},p}$ is
totally ordered by the relation $\geq$.

\item[-] ${\cal F}$ is not linearizable at $p$. Then there exist
local analytic coordinates $u,v$ providing the Poincaré-Dulac normal
form of ${\cal F}_p$, which is $(nu+v^n)du-vdv$ for a certain
positive integer $n$. Using this fact, it is straightforward to see
that, as in the previous case, the set ${\cal C}_{{\cal F},p}$ is
totally ordered by the relation $\geq$. Observe that this case does
not occur when ${\cal F}$ admits a rational first integral.
\end{itemize}

Each set ${\cal C}_{{\cal F},p}$ defines, then, a valuation of the
fraction field of $\co_{\gp^2,p}$ and a simple complete primary
ideal $I_p$ of that local ring \cite[pages 389 to 391]{zar}; in
fact, when ${\cal F}$ is linearizable at $p$,
$I_p\co_{\gp^2,p}^{an}$ coincides with the integral closure of
${J}_p$ in $\co_{\gp^2,p}^{an}$. Set $D_p({\cal F})$ the unique
exceptionally supported divisor on $X_{\cal F}$ such that
$I_p\co_{X_{\cal F}}=\co_{X_{\cal F}}(-D_p({\cal F}))$ and, for each
positive integer $d$ and each map $\bold{k}: NRed(\cf)\rightarrow
\gz_+$, denote by ${\cal L}({\cal F},d,\bold{k})$ the sheaf
${\pi_{\cal F}}^*\co_{\gp^2}(d)\otimes \co_{X_{\cal F}}(-\sum
\bold{k}(p) D_p({\cal F}))$, where the sum is taken over the set of
non-reduced singularities of ${\cal F}$.

We shall use the above notations in the statement and proof of the
next result, that characterizes those non-degenerated plane
foliations of degree $r\not=1$ admitting a rational first integral.

\begin{teo}\label{principal}
Let ${\cal F}$ be a non-degenerated foliation on $\gp^2$ of degree
$r\not=1$ and let $\{(\delta_p , \rho_p)\}_{p\in NRed({\cal F})}$ be
the set of pairs of eigenvalues associated with the non-reduced
singularities. Then, ${\cal F}$ has a rational first integral if and
only if there exist a positive integer $d$ and a map $\bold{k}:
NRed(\cf)\rightarrow \gz_+$ such that:

\begin{itemize}
\item[(a) ] $d^2=\sum_{p\in NRed({\cal F})} \bold{k}(p)^2 \rho_p
\delta_p$.

\item[(b) ] $d(r+2)=\sum_{p\in NRed({\cal F})}
\bold{k}(p)(\rho_p+\delta_p)$.

\item[(c) ]  $h^0(X_{\cal F},{\cal L}({\cal F},d,\bold{k}))=2$.

\item[(d) ] There exist two curves of ${\pi_{\cal F}}_* |{\cal
L}({\cal F},d,\bold{k})|$ (say, with equations $H_1=0$ and $H_2=0$)
such that, for each $p\in NRed({\cal F})$,  the germs at $p$ of the
curves with equations $H_1=0$ and $H_2=0$ are equisingular and the
one of $H_1H_2=0$ satisfies the following conditions: it is reduced,
its associated Milnor and Tjurina numbers coincide and its
topological singularity type is $S(\rho_p,\delta_p,2\bold{k}(p))$.

\item[(e) ] Each point in $Red({\cal F})$ is a singular point of
some curve in ${\pi_{\cal F}}_*|{\cal L}({\cal F},d,\bold{k})|$.

\end{itemize}

Moreover, in this case, a primitive rational first integral of
${\cal F}$ is $F/G$, where $F$ and $G$ are two homogeneous
polynomials defining any two different curves in ${\pi_{\cal
F}}_*|{\cal L}({\cal F},d,\bold{k})|$.

\end{teo}

\noindent {\it Proof}. Assume first that ${\cal F}$ has a rational
first integral and let $R=F/G$ be a primitive one. Set $d$ the
degree of $F$ and $G$. The morphism $\pi_{\cal F}$ above defined is
the one eliminating the base points of the pencil ${\cal P}(R)$ or,
equivalently, eliminating the indeterminacies of the rational map
$\gp^2 \cdots \rightarrow \gp^1$ provided by $R$ \cite[II.7]{beau}.
Moreover, for each $p\in NRed({\cal F})$, if $k_p$ denotes the
number of branches through $p$ of a general curve $C$ of the pencil
${\cal P}(R)$, the multiplicities $m_q$ of the strict transforms of
$C$ at the points $q\in {\cal C}_{{\cal F},p}$ coincide with the
ones of the strict transforms of any general element of the ideal
$I_p^{k_p}$ (since $I_p\co_{\gp^2,p}^{an}$ is the integral closure
of ${J}_p$). Thus, taking the function $\bold{k}: NRed({\cal
F})\rightarrow \gz_+$ such that $\bold{k}(p)=k_p$ for all $p\in
NRed({\cal F})$, the following equality of sheaves holds: ${\cal
L}({\cal F},d,\bold{k})=\pi_{\cal F}^* \co_{\gp^2}(d)\otimes
\co_{X_{\cal F}}(\sum m_q E_q^*)$, where the sum is taken over the
points $q\in {\cal C}_{\cal F}$ and $E_q^*$ denotes the pull-back on
$X_{\cal F}$ of the exceptional divisor appearing in the blowing-up
centered at $q$. Since ${\cal P}(R)$ coincides with the pencil
${\cal P}_{\cal F}$ given in \cite[Lem. 1]{gal-mon}, we can apply
this result to deduce that ${\cal P}(R) = {\pi_{\cal F}}_*|{\cal
L}({\cal F},d,\bold{k})|$; hence $(c)$ is satisfied. Clauses $(d)$
and $(e)$ follow from Proposition \ref{primero}, $(a)$ is a
consequence of applying Bézout's Theorem to two general curves in
${\cal P}(R)$ and $(b)$ follows from a result in \cite{poi2}.

Conversely, assume the existence of a positive integer $d$, a map
$\bold{k}: NRed(\cf)\rightarrow \gz_+$ and curves $H_1=0$ and
$H_2=0$ as in the statement. Set $R=H_1/H_2$ and let ${\cal H}$ be
the foliation provided by the derivation of the rational function
$R$. Using Condition $(d)$ and similar arguments to those given in
the last part of the proof of Proposition \ref{primero} one can
deduce that, for each $p\in NRed({\cal F})$, the Milnor number
$\mu_p({\cal H})$ equals $1$, the germs at $p$ of $H_1=0$ and
$H_2=0$ have singularity type $S(\rho_p,\delta_p,\bold{k}(p))$ and,
moreover, they have the same minimal resolution of singularities.
From the last assertion one has also that
$i_p(H_1=0,H_2=0)=\bold{k}(p)^2\rho_p\delta_p$, where $i_p$ stands
for the intersection multiplicity at $p$. This fact, Bézout's
Theorem and Condition $(a)$ show that the base points on $\gp^2$ of
the pencil $\cp(R)={\pi_{\cal F}}_*|{\cal L}({\cal F},d,\bold{k})|$
(the equality holds from $(c)$) are exactly those in $NRed({\cal
F})$. Hence, $NRed({\cal H})=NRed({\cal F})$.

It follows from the above paragraph that the pair of eigenvalues
attached to each non-reduced singularity $p$ of ${\cal H}$ is
$(\delta_p,\rho_p)$ and the number of branches through $p$ of a
general curve of the pencil ${\cal P}(R)$ is $\bold{k}(p)$. Then, if
$r'$ denotes the degree of the foliation ${\cal H}$, the equality
 $d(r'+2)=\sum_{p\in NRed({\cal F})} \bold{k}(p)(\rho_p+\delta_p)$ holds
\cite{poi2}. By Condition $(b)$, we deduce that $r=r'$. Therefore,
\begin{equation}\label{milnor2}
\sum_{p\in Sing({\cal H})} \mu_p({\cal H})=r^2+r+1
\end{equation}
and so, $\# Red({\cal H}) \leq \# Red({\cal F})$, where $\#$ stands
for the cardinality. But, by Condition $(e)$ and Part $(b)$ of Lemma
\ref{lema}, each point in $Red({\cal F})$ is also a singular point
of ${\cal H}$. Hence one gets that $Ind(R)=NRed({\cal H})=NRed({\cal
F})$, $Sing(R)=Red({\cal H})=Red({\cal F})$ and, taking into account
(\ref{milnor2}), ${\cal H}$ is a non-degenerated foliation. Thus, by
Proposition \ref{primero}, $R$ is a rational first integral of
${\cal F}$.

Finally, since $h^0(X,{\cal L}({\cal F},d,\bold{k}))=2$ we deduce
that
 $R$ is a primitive rational first integral of ${\cal F}$
 by the proof of \cite[Th. 2 (a)]{gal-mon}.\findemo

\begin{rem}
{\rm Note that, if Condition $(d)$ in the statement of the above
result is satisfied for two curves of ${\pi_{\cal F}}_*|{\cal
L}({\cal F},d,\bold{k})|$, then it holds for any pair of general
elements of this linear system. }
\end{rem}

\begin{exa}
{\rm The pair of Diophantine equations of a fixed non-degenerated
foliation $\cal F$ given by  clauses (a) and (b) of Theorem
\ref{principal}, where the unknowns are $d$ and $ \{\bold{k}(p)\}_{p
\in NRed({\cal F})}$, can have infinitely many solutions. Indeed,
foliations in the set $\{{\cal F}^r_\alpha\}_{\alpha \in E}$, $2
\leq r \leq 4$, given in \cite{l-n}, ($E$ is a countable and dense
set of parameters in the set of complex numbers) have singularities
of fixed analytic type \cite[Def. 1]{l-n},  rational first integral
and their degrees can be chosen arbitrarily large. Therefore we get
Diophantine equations
\[
\begin{array}{c||c||c}
  d^2= 2 \sum_{i=1}^2 k_i^2 +6 \sum_{i=3}^5 k_i^2 \;\;&\;\;
  d^2=  \sum_{i=1}^3 k_i^2 + 2 \sum_{i=4}^8 k_i^2 \;\;&\;\; d^2=  \sum_{i=1}^{12}k_i^2 \\
  4d = 3 \sum_{i=1}^2 k_i +5 \sum_{i=3}^5 k_i \;\;&\;\; 5d = 2 \sum_{i=1}^3 k_i +3 \sum_{i=4}^8 k_i
  \;\;& \;\;
  3d =  \sum_{i=1}^{12}k_i
\end{array}
\]
associated, respectively, to the cases $r=2,3,4$, where for
simplicity we have set $\bold{k}(p_i) = k_i$, that have infinitely
many solutions (attached to the corresponding first integrals).

To end this example, we examine clauses in Theorem \ref{principal}
for the foliation ${\cal F}^4_0$ that is given by the 1-form $A dX
+B dY + C dZ$, where $A=(Y^3-Z^3)YZ$, $B=(Z^3-X^3)XZ$ and
$C=(X^3-Y^3)XY$. It has 12 points in $NRed({\cal F})$ and 9 points
in $Red({\cal F})$. As ${\cal F}^4_0$ has rational first integral,
clauses (a) to (e) in Theorem \ref{principal} must hold. In fact,
$d=6$, $k_1 = k_{2} = k_{3} =3$ and  $k_i = 1$, $4 \leq i \leq 12$
are solutions for the above last Diophantine equations. Forms $H_1=
3 X^3Y^3- X^3Z^3- 2 Y^3Z^3$ and $H_2= 2X^3Y^3- X^3Z^3-  Y^3Z^3$ span
the vector space $H^0(X_{\cal F},{\cal L}({\cal F},6,\bold{k}))$ for
$\bold{k}$ defined as above and a suitable ordering of the points in
$NRed({\cal F})$. The equation of the germ of curve defined by $H_1$
($H_2$, respectively) at the point $p=(0:0:1)\in NRed(\cf)$ (whose
image by $\bold{k}$ equals 3) in suitable local coordinates is given
by $h_1=x^3+2y^3-3x^3y^3$ ($h_2=x^3+y^3-2x^3y^3$, respectively). The
germ defined by $h_1h_2$ factorizes into a product of 6 smooth
transversal analytically irreducible germs passing through $p$
(having, thus, singularity type $S(1,1,6)$) and its Milnor and
Tjurina numbers are both equal to 25 (we have used {\sc Singular}
\cite{sin} to do the computation); obviously the germs defined by
$h_1$ and $h_2$ are equisingular. The same situation happens for the
two remaining points of $NRed(\cf)$ whose image by $\bold{k}$ is 3.
The germs of $H_1$ and $H_2$ at any point $q\in NRed(\cf)$ such that
$\bold{k}(q)=1$ are analytically irreducible, smooth and
transversal. This shows that the conditions given in (d) are
satisfied for $H_1$ and $H_2$. Finally  (e) also holds because the
points $(j^r:0:1)$ ($(1:j^r:0)$, respectively) ($(0:j^r:1)$,
respectively), where $r \in \{0,1,2\}$ and $j:=e^{2 \pi i /3}$, are
singular points of $Y^3(X^3-Z^3)$ ($Z^3(Y^3-X^3)$, respectively)
($X^3(Y^3-Z^3)$, respectively).}
\end{exa}

\begin{rem} \label{nota2}
{\rm The unique non-degenerated foliations of degree one are
defined, up to projective isomorphism, by the differential 1-forms
$aYZ\;dX+bXZ\;dY - (a+b)XY\;dZ$, where $a$ and $b$ are positive
integers, $(X:Y:Z)$ being projective coordinates on the complex
projective plane (see \cite[Sect. 4]{cam-oli}). All of them are
algebraically integrable, since $X^a Y^b /Z^{a+b}$ is a rational
first integral.

}
\end{rem}

\begin{rem}
{\rm  Theorem \ref{principal} allows to decide whether a
non-degenerated foliation of degree $r \neq 1$ $\cal F$, defined by
a projective differential 1-form $\Omega$, has or not a rational
first integral $\frac{F}{G}$ such that the degree of $F$ and $G$ is
less than a fixed value $t$ (and to compute it, if it exists). The
procedure should be as follows: firstly, perform a minimal
resolution of ${\cal F}$ and compute the divisors $D_p({\cal F})$
for all $p\in NRed({\cal F})$. Secondly, consider the integers $d<t$
and the finite set of maps $\bold{k}:NRed(\cf)\rightarrow \gz_+$
satisfying the conditions $(a)$, $(b)$ and $(c)$ of Theorem
\ref{principal}. Then, ${\cal F}$ has a rational first integral with
the above condition if and only if $d(\frac{F}{G})\wedge \Omega=0$
for some basis $\{F,G\}$ of a linear system of the type ${\pi_{\cal
F}}_*|{\cal L}({\cal F},d,\bold{k})|$.

}
\end{rem}

\begin{exa}
{\rm Take projective coordinates $(X:Y:Z)$ on the complex projective
plane and consider the non-degenerated foliation $\cal F$ defined by
the projective differential 1-form $\Omega=AdX+BdY+CdZ$, where
$$A=Z(2X^3+2Y^3-YZ^2), \;\;\;
B=Z(XZ^2-6XY^2)\;\;\; \mbox{and}\;\;\; C=4XY^3-2X^4.$$ It has 6
non-reduced singularities $p_1,p_2,\ldots,p_6$ and the morphism
$\pi_{\cal F}:X_{\cal F}\rightarrow \gp^2$ is the composition of the
blowing-ups with centers $p_1, p_1',p_2,p_2',p_3,p_3',p_4,p_5$ and
$p_6$, where $p_i'$ is infinitely near to $p_i$ for $1\leq i\leq 3$.
The  pair of eigenvalues associated with $p_i$, $1\leq i\leq 3$
($4\leq i\leq 6$, respectively), is $(1,2)$ ((1,1), respectively).
Using results of \cite{lipman} we see that the divisor
$D_{p_i}({\cal F})$  is $E^*_{p_i}+E^*_{p_i'}$ ($E^*_{p_i}$,
respectively) for $1\leq i\leq 3$ ($4\leq i \leq 6$, respectively),
where $E^*_q$ denotes the total transform on $X_{\cf}$ of the
exceptional divisor appearing in the blowing-up centered at $q$. The
integer $d=3$ and the map $\bold{k}\equiv 1$  satisfy the conditions
$(a)$, $(b)$ and $(c)$ of Theorem \ref{principal} and the linear
system ${\pi_{\cf}}_* |{\cal L}(\cf,3,\bold{k})|$ is spanned by the
curves with equations $F:=X^3-2Y^3+YZ^2=0$ and $G:=XZ^2=0$. After
checking it, we assert that $F/G$ is a rational first integral of
$\cf$.
 }
\end{exa}

When a non-degenerated plane foliation ${\cal F}$ has a rational
first integral, the cardinality of the set $NRed({\cal F})$ of
non-reduced singularities is not arbitrary, as the next result
shows. Notice that this gives an easy to check criterion to decide
that a foliation has no rational first integral.

\begin{teo} \label{cota}
Let ${\cal F}$ a degree $r$ non-degenerated foliation on $\Pp^2$
which has a rational first integral. Then
$$
r+1 \leq n,
$$
where $n$ is the number of non-reduced singularities of ${\cal F}$.

\end{teo}

\noindent {\it Proof.} Let $\{(\delta_p , \rho_p)\}_{p\in NRed({\cal
F})}$ be the set of pairs of eigenvalues associated with the
non-reduced singularities of the foliation ${\cal F}$. Consider a
primitive rational first integral of ${\cal F}$, $R$, and let $C$ be
a general curve of the pencil $\cp(R)$. For each $p\in NRed(\cf)$,
let $k_p$ be the number of branches of the singularity of $C$ at
$p$. Assume that $C$ has degree $d$. Bertini's Theorem proves that
the unique singularities of $C$ are the non-reduced singularities of
${\cal F}$ and similar arguments to some used in the proof of
Proposition \ref{primero} show that each non-reduced singularity $p$
is quasi-homogeneous with associated weights $(k_p\rho_p,
k_p\delta_p)$. So, in suitable coordinates $\{x, y\}$ in ${\cal
O}_{\Pp^{2},p}^{an}$, the germ of $C$ at $p$ is given by
$H(x^{\rho_p},y^{\delta_p})$, $H$ being an homogeneous polynomial in
two indeterminates of degree $k_p$. So the Milnor and Tjurina
numbers coincide at each singularity of $C$.

Set $\mu$ the sum of Milnor numbers of the singularities of $C$.
Then,
\begin{equation}
\label{milnor} (d-1)(d-r-1) \leq \mu
\end{equation}
by a result of du Plessis and Wall \cite{p-w}. $C$ is a general
curve of the pencil ${\cal P }(R)$ and, so, the Milnor-Jung formula
states
\[
\mu_{p} = 2 {\bar \delta}_{p} - k_p +1
\]
for all $p\in NRed({\cal F})$, where $\mu_{p}=\mu(C,p)$ (${\bar
\delta}_{p}$, respectively) is the Milnor number (genus diminution
or $\delta$-invariant, respectively) of $C$ at the singularity
$p$. Hence,
\[
\mu = \sum_{p\in NRed({\cal F})}(2 {\bar \delta}_{p} - k_p) +n.
\]
On the one hand, by the genus formulae, one gets:
\[
2g -2 = C \cdot C + C \cdot K_{\Pp^2} - \sum_{p\in NRed({\cal F})}
k_p^2 \rho_p \delta_p + \sum_{p\in NRed({\cal F})} k_p (\rho_p +
\delta_p -1),
\]
where $g$ is the geometrical genus of $C$ and $K_{\Pp^2}$  a
canonical divisor of the projective plane (see \cite[pages 279 and
280]{g-h}). On the other hand, it holds
\[
2g = (d-1)(d-2) - 2 \sum_{p\in NRed({\cal F})} {\bar \delta}_{p}
\]
(see, for instance, \cite[Sect. IV, Ex. 1.8]{har}). Then, we obtain
the following equality:
\[
\mu = \sum_{p\in NRed({\cal F})} [k_p^2 \rho_p \delta_p -
 k_p (\rho_p + \delta_p)] + n.
\]
Finally, $$\mu = d^2 -(r+2)d +n,$$ since $d^2=\sum_{p\in NRed(\cf)}
k_p^2 \rho_p \delta_p$ and $\sum_{p\in NRed(\cf)} k_p (\rho_p +
\delta_p)=(r+2)d$ \cite{poi2}, which concludes the proof after
replacing in
(\ref{milnor}) the value of $\mu$. \findemo\\

\begin{exa}
{\rm For each $a\in \gc\setminus \{0,1,-1\}$ consider the
non-degenerated plane foliation $\cf_a$ defined by the projective
1-form $\Omega=AdX+BdY+CdZ$, $(X:Y:Z)$ being projective coordinates
in the complex projective plane, where:
$$A=Z(aXZ-Y^2+Z^2), \;\;\;
B=Z(X^2-Z^2)\;\;\; \mbox{and}$$ $$C=XY^2-aX^2Z-XZ^2-X^2Y+YZ^2.$$ The
cardinality of $NRed({\cal F}_a)$ is $1$ or $2$, depending on the
value of the parameter $a$ (see \cite[Example 5]{gal-mon}). Since
the degree of the foliation is 2, Theorem \ref{cota} allows to
discard the existence of a rational first integral for any foliation
in  the family $\{{\cal F}_a\}_{a\in \gc\setminus \{0,1,-1\}}$.

}
\end{exa}

Finally, we give an example which shows that the lower bound of
Theorem \ref{cota} turns out to be sharp.

\begin{exa}
{\rm Consider any foliation ${\cal F}$ of the family given in Remark
\ref{nota2}. ${\cal F}$ has degree 1 and it has three singular
points: $(0:0:1), (0:1:0)$ and $(1:0:0)$. A rational first integral
of ${\cal F}$ is $\frac{X^a Y^b}{Z^{a+b}}$ (for suitable $a$ and
$b$) and, so, the non-reduced singularities are $(0:1:0)$ and
$(1:0:0)$.

}
\end{exa}

\vspace{3mm}
\par
\noindent \footnotesize \textsc{C. Galindo:}  Universidad Jaume I,
Dep. de Matemáticas (ESTCE), Campus Riu Sec, 12071  Castellón, Spain. \\
E-mail address: \verb"galindo@mat.uji.es"\\[0.75mm]

\noindent \textsc{F. Monserrat:}  Universidad Politécnica de
Valencia, Instituto Universitario de Matemática Pura y Aplicada,
Camino de Vera s/n, 46022 Valencia, Spain. \\
E-mail address: \verb"framonde@mat.upv.es"

\end{document}